\documentclass[12pt]{amsart}

\makeatletter                                            
\renewcommand{\@begintheorem}[2]{                        
\rm \trivlist \item [\hskip \labelsep {\bf #2\ \ #1.}]   
                                }                        
\makeatother                                             

\makeatletter

\newcommand{\ts}{\vspace{\baselineskip}\noindent{\bf Proof.}$\;\;$}
\newcommand{\ZZ}{{\bf Z}}
\newcommand{\QQ}{{\bf Q}}

\newcommand{\CC}{{\bf C}}

\newcommand{\FF}{{\bf F}}

\newcommand{\PP}{{\bf P}}

\newcommand{\Km}{\mbox{Km}}
\newcommand{\et}{\text{\'et}}
\newcommand{\tpi}{\mbox{$\frac{1}{t+1}$}}
\newcommand{\dd}{\text{d}}

\begin{document}

\title{An isogeny of $K3$ surfaces}
\author{Bert van Geemen and Jaap Top}
\address{Dipartimento di Matematica, Universit\`a di Milano,
  via Saldini 50, I-20133 Milano, Italia}
\address{IWI, Rijksuniversiteit Groningen, P.O.Box 800, 9700 AV Groningen,
the
Netherlands}
\email{geemen@mat.unimi.it}\email{top@math.rug.nl}

\begin{abstract}
In a recent paper Ahlgren, Ono and Penniston described the L-series
of $K3$ surfaces from a certain one parameter family in
terms of those of a particular family of elliptic curves.
The Tate conjecture predicts the existence of a correspondence
between these $K3$ surfaces and certain Kummer surfaces related to these
elliptic curves.
A geometric construction of this correspondence is given here,
using results of D.~Morrison on Nikulin involutions.
\end{abstract}

\maketitle

\section{The family}

\subsection{}\label{1.1}
Recently, Ahlgren, Ono and Penniston \cite{AOP} studied
the $K3$ surfaces $X_t$ which are the minimal resolutions of double
covers of
$\PP^2$
branched over a union of $6$ lines (hence over a sextic curve):
$$
X_t:\quad y^2=xz(x+1)(z+1)(x+zt).
$$
Using an elaborate but elementary calculation with
character sums, they determined the zeta function
of $X_t/\FF_p$. One way of interpreting their
result is as follows.

For general $t\in\QQ$ the N\'eron-Severi group of $X_t$ has rank $19$
(cf.\ Lemma~\ref{Tgen} below).
Hence there is an isomorphism of
$G_\QQ={\rm Gal}(\overline{\QQ}/\QQ)$ representations:
$$
H^2_{\et}(X_{t,\overline{\QQ}},\QQ_\ell)\cong
T_{t,\ell}\oplus \QQ_\ell(-1)^{19}
$$
for some $\ell$-adic representation $T_{t,\ell}$ of dimension $3$.

Consider the elliptic curve $E_t$ and its quadratic twist $E_t^{(t+1)}$:
\[
E_t:        \quad y^2=(x-1)(x^2-\tpi),\qquad
E_t^{(t+1)}:\quad (t+1)y^2=(x-1)(x^2-\tpi).
\]
The Kummer surface $\Km(E_t\times E_t^{(t+1)})$
is by definition the smooth surface
obtained by blowing up the $16$
double points of the quotient $(E_t\times E_t^{(t+1)})/([-1]\times[-1])$.
This Kummer surface is also a $K3$ surface.
Since $E_t^{(t+1)}$ is a quadratic twist of $E_t$,  we obtain another
3-dimensional $G_\QQ$-representation
$$
\mbox{Sym}^2\left(H^1_{\et}(E_{t,\overline{\QQ}},\QQ_{\ell})\right)
(\chi^{(t+1)})
\qquad(\subset
H^2_{\et}(\Km(E_t\times E_t^{(t+1)})_{\overline{\QQ}},\QQ_\ell)).
$$
Here $\chi^{(t+1)}$ is the Dirichlet character of the quadratic
extension $\QQ(\sqrt{t+1})/\QQ$ if $t+1$ is not a square in $\QQ$,
else it is trivial.

\vspace{\baselineskip}\noindent
{\bf Proposition (Ahlgren, Ono, Penniston)}.
{\sl With notations as above, the two Galois representations
$T_{t,\ell}$ and
$\mbox{Sym}^2\left(H^1_{\et}(E_{t,\overline{\QQ}},\QQ_{\ell})\right)
(\chi^{(t+1)})$
are isomorphic.}

\vspace{\baselineskip}

This isomorphism produces, via the K\"unneth formula and Poincar\'e duality,
a Galois invariant class in
$H^4_{\et}(X_{t,\overline{\QQ}}\times
\Km(E_t\times E_t^{(t+1)})_{\overline{\QQ}},\QQ_\ell)$.
The Tate conjecture asserts that for
a variety $Z$, defined over $\QQ$, the subspace of
Galois invariants in
$H^{2p}_{\et}(Z_{\overline{\QQ}},\QQ_\ell)(p)$
is spanned by classes of codimension $p$ cycles defined over $\QQ$.
Combined with the proposition above this suggested our main result:

\subsection{Theorem}\label{main}
{\sl
For $t\in\QQ$ there exists a
correspondence
$$
\Gamma_t\;\subset\; X_{t}\times \Km(E_t\times E_t^{(t+1)}),
$$
defined over $\QQ$, which induces an isomorphism of $G_\QQ$-representations:
$$
[\Gamma_t]:
T_{t,\ell}\;\stackrel{\cong}{\longrightarrow}\;
\mbox{Sym}^2\left(H^1_{\et}(E_{t,\overline{\QQ}},\QQ_{\ell})\right)
(\chi^{(t+1)}).
$$
}

\vspace{\baselineskip}\noindent
In Remark 4.4 of the paper \cite{AOP}, the authors suggest finding a
dominant rational map from $X_t$ to $K_t=Km(E_t\times E_t)$.
This is actually possible, but only over a finite extension of $\QQ(t)$,
and we do produce such a map.

\subsection{Proposition}\label{prp}
{\sl Let $K$ be a field of characteristic $\neq 2$ and take $t\neq 0,-1$ in
$K$.
Then there exists a dominant rational map from $X_t$ to
$K_t=Km(E_t\times E_t)$ over a finite extension of $K$.
}

\subsection{}
Such a geometric relation (at least, over the
complex numbers) between the two
families of $K3$ surfaces can also be shown to exist using their Picard-Fuchs
differential equations. This has been worked out by Ling Long \cite{L}.

We now briefly outline the general facts we used and the strategy we
followed to obtain our result.

\subsection{General results.}\label{general}
The general $X_t$ has a N\'eron-Severi group $NS(X_t)$ of rank~$19$, and thus
its transcendental lattice $T=NS(X_t)^\perp$ has rank~$3$,
we computed (cf.\ Lemma~\ref{Tgen})
that
$$
T\cong \langle 2\rangle \oplus \langle 2\rangle \oplus \langle -2\rangle.
$$

Recall that $\Km(A)$, the Kummer surface of an abelian surface $A$,
is the $K3$ surface obtained by blowing up  the $16$ singular points
of the quotient of $A$ by the involution $a\mapsto -a$.
A $K3$ surface $S$ with rank $NS(S)=19$ is a Kummer surface if and only if
the (even) quadratic form $q:T_S=NS(S)^\perp\rightarrow 2\ZZ$
obtained from the
intersection product on $H^2(S,\ZZ)$ has values in $4\ZZ$
(cf.\ \cite[Prop.~4.3]{morrison}).
In particular, the general $X_t$ is not a Kummer surface.

The transcendental lattice $T_S$ of any $K3$ surface $S$ of rank $19$
embeds into $U^3$ (\cite[Cor.~2.6]{morrison},
where $U$ is the hyperbolic plane ($\ZZ^2$ with quadratic form
$q(x)=2x_1x_2$).
This gives an embedding of $T_S$ in the $K3$ lattice $U^3\oplus E_8(-1)^2$
which is unique up to isometry (\cite[Cor.~2.10]{morrison}).
The N\'eron-Severi group $NS(X)$ of $X$ thus contains $E_8(-1)^2$.
Theorem 5.7 of \cite{morrison} now implies that $X$ has a Nikulin involution
$\iota$ (that is, an involution which acts trivially on $H^{2,0}(X)$).
The involution has $8$ fixed points, blowing them up and taking the quotient
we get a $K3$ surface $V$ with $T_V\cong T(2)$ (\cite[Thm.~5.7ii]{morrison}).
Hence $V$ is a Kummer surface.
The corresponding abelian surface $A$ has transcendental lattice $T_A\cong T$
(\cite[Prop.~4.3]{morrison}). The following diagram summarizes the situation,
it is called a Shioda-Inose structure for $X$.
\[\begin{array}{ccccc}
X &&&& A \\
 & \searrow && \swarrow & \\
&& X/\iota \cong \Km(A) &&
\end{array}
\]

\subsection{Summary.}
For the general $X_t$ it is rather easy to find a sublattice $E_8(-1)^2$
of $NS(X_t)$, see section \ref{e8s}. Since the Nikulin involution exchanges
the two copies of $E_8(-1)$ and we have an interpretation of the simple
roots as nodal curves on the surface, we can make an educated guess as to
what the involution should be. In section \ref{2.5} we give the involution
explicitly and we determine the quotient $K3$ surface $V_t$.
In section \ref{vkum}, we show that $V_t$ is isomorphic to a double cover
$W_t$ of $\PP^2$ branched along 6 lines which are tangent to a conic.
This shows that $V_t\cong \Km(JC_t)$ where $JC_t$
is the Jacobian of the genus two curve $C_t$
which is the double cover of the conic
branched in the 6 points of tangency of the lines to the conic.
The abelian surface $JC_t$ is isogenous to a product of two elliptic
curves $F_t\times F_t'$ which are quadratic twists of $E_t$
(section \ref{4.5}).
A main problem is that most isogenies and isomorphisms are not defined
over $\QQ$ (or $\QQ(t)$).
The varieties involved do have models over $\QQ$, but one has to choose the
right one (or twist a given one) so as to have a non-trivial
correspondence defined over $\QQ$.
We conclude with some observations on the `famous'
$K3$ surface $X_{-1}$.

\subsection{Previous work.}
In the literature, several results comparable to Proposition~\ref{prp}
can be found. However, we are not aware of any cases except the present one
where an explanation
is given how such isogenies may be constructed. We mention some examples
here.
Note
that they are older than Morrison's paper which provides the basic technique
for our construction. It may be interesting to study whether Long's method
mentioned above can be used in the following examples as well to predict the
existence of
the isogenies involved.

In 1977, M.~Mizukami \cite{Mi} showed that the Kummer surface
$\Km(E'_t\times E'_t)$ is isogenous to the $K3$ surface $X'_t$,
for $t\neq \pm1$ in $\CC$, where
\[X'_t\;:\quad x_1^4+x_2^4+x_3^4+x_4^4+2t(x_1^2x_2^2+x_3^2x_4^2)=0
\]
and
\[E'_t\;:\quad y^2=(x^2+1)(x^2+(1-t)/2).
\]
This is proven by explicitly giving a rational $4:1$ map from
$E'_t\times E'_t$ to $X'_t$.

Similarly, W.~Hoyt \cite{Hoyt} in 1984 presents an explicit rational
dominant map from the product $E''_t\times E''_t$, where
\[E''_t\;:\quad y^2=x^3-(12-9t)x+16-18t,\]
to the $K3$ surface $X''_t$ corresponding to the equation
\[ s^2=x(x-1)(x-t)y(y-1)(y-x).\]

\section{The $K3$ surfaces $X_t$}

\subsection{Singularities of the branch curve.}\label{singbra}
For $t\neq 0$, the branch curve of the double cover defining $X_t$
consists of 6 lines (including the line at infinity), from now on we
assume $t\neq 0$.

For $t\neq -1$ these lines meet in
6 double points and 3 triple points.
To obtain the corresponding $K3$ surface, one blows up the double
and triple points. Over a triple point, one must next blow up the three
intersection points of the strict transforms of the three lines
and the exceptional divisor. We denote by $E_P$
the inverse images in the $K3$ surface
$X_t$ of the strict transform of the fibre of the first blow up in $P$.
Furthermore, by $E_P^{l=0}$ we denote the inverse image of the exceptional
divisor over the point of intersection of $E_P$ and the strict transform of
the line $l=0$. For symmetry reasons we write the point with coordinates
$(x:z:1)\in \PP^2$ as $(x:z:-1)$,
thus the exceptional divisor over the double
point $(x,z)=(-1,0)$ is denoted by $E_{101}$.
All these curves,
as well as the inverse images of the lines which
make up the sextic (we denote these simply by $l=0$ as in $\PP^2$),
are smooth rational curves, hence $(-2)$-curves, in the $K3$ surface.
See \cite{AOP}, p.\ 363, figure 1, for a picture of
the intersection graph of these $-2$-curves.

In the special case $t=-1$ there are 3 double points and 4 triple points.
The $6+3+4\cdot 4=25$ rational curves in $X_{-1}$ with self intersection $-2$
obtained as above are denoted in the same way.

\subsection{The case $t=-1$.}\label{x1}
It is shown in \cite[p.~298]{persson} (see also the proof of \ref{Tgen})
that the $K3$ surface $X_{-1}$ has transcendental lattice $T_{-1}$
of rank~$2$  and discriminant~$4$, hence $T_{-1}$ must be
$$
T_{-1}= \langle 2\rangle\oplus \langle 2\rangle.
$$

Vinberg \cite{vinberg} studied the (unique) $K3$ surface with this
transcendental lattice and observed (\cite[2.1]{vinberg})  that its Picard
lattice is isomorphic to the
sublattice of $\ZZ^{20}$ (with quadratic form $x_1^2-\sum_{i=2}^{20}x_i^2$)
given by the vectors $x$ with $\sum x_i\equiv 0\pmod 2$.
Shioda and Inose \cite{SI}
showed that $X_{-1}$ is the desingularisation of the quotient of
$E_i^2$, the self product of the elliptic curve $E_i=\CC/\ZZ[i]$,
by the automorphism $\phi$ of order $4$ induced by
$(z_1,z_2)\mapsto (iz_1,-iz_2)$ on $\CC^2$, see also Section~\ref{a1}.

The following lemma is not used in the proof of the main result,
but it does show that $X_t$ is not a Kummer surface,
hence we cannot avoid the Nikulin involution.

\subsection{Lemma.}\label{Tgen}
{\sl The N\'eron-Severi group of the general $X_t$ has rank $19$ and it
is generated by nodal curves defined over $\QQ(t)$.
The transcendental lattice $T$ of the general $X_t$ is given by
$$
T\cong \langle 2\rangle \oplus \langle 2\rangle \oplus \langle -2\rangle.
$$
}

\ts
First we consider the special case $t=-1$. The sublattice of
$H^2(X_{-1},\ZZ)$ generated by the following $20$ nodal curves:
$E_{011}, E_{010}^{x=0}, E_{111}, E_{111}^{x=-1}$ and the 16 curves
which span the two copies of $E_8$ given in section \ref{e8s} has
rank 20 and determinant~$-4$, as can be verified by a computation
with their
intersection matrix. Hence the determinant of the transcendental lattice
$T_{-1}$ is either $4$ or $1$, but since $T$ is positive definite and
even, $det(T_{-1})=4$ (and in fact 
$T_{-1}\cong \langle2\rangle\oplus \langle 2\rangle$).
It follows that the $20$ curves are a $\ZZ$-basis of $NS(X_{-1})$.
According to Nikulin, \cite[Thm.~1.14.4]{Ni}, the embedding
of $T$ in the $K3$ lattice $U^3\oplus E_8(-1)^2$ is unique up to
isometry. We fix such an embedding and identify $NS(X_{-1})$ with
$T^{\perp}$.

Since the $K3$ surfaces $X_t$ depend on one parameter, $NS(X_t)$ has
rank at most 19 for general $t$.
For $t\neq -1$, there are three points $p_a=(-1,t^{-1}),p_b=(t,-1),
p_c=(-1,-1)$
which specialize to the triple point $(x,z)=(-1,-1)$ as $t\rightarrow -1$.
It is easy to check that the lattice spanned by the
24 nodal curves in $X_t$ obtained from
the desingularization of the branch locus is isomorphic to the sublattice
of $NS(X_{-1})$ spanned by the corresponding curves in case $t=-1$
and where $E_a$ maps to $E_{111}^{x=-1}+E_{111}+E_{111}^{x=-zt}$ et cetera
(so the exceptional curve over the intersection of two lines $l=0$ and $m=0$
maps to the sum of the curves $E_{111}^{l=0}+E_{111}+E_{111}^{m=0}$).
Now it is easy to see that the image of $NS(X_t)$ in $NS(X_{-1})$ is
contained in the orthogonal complement $L$ of $E_{111}$. Hence $L$ is
a primitive sublattice and it is actually generated by the $24$ nodal curves
(in fact, the 16 curves in the copies of $E_8$ and
$E_{011}, E_{010}^{x=0}, E_c$ are a $\ZZ$-basis of $L$).
Therefore $NS(X_t)\cong L$ for general $t$. A computation shows that
$L\cong E_8(-1)^2\oplus \langle-2\rangle^2\oplus \langle2\rangle$ and
that $L^\perp$ is isomorphic to 
$\langle 2\rangle \oplus \langle 2\rangle \oplus \langle -2\rangle $.
\hfill $\Box$

\section{The isogeny}

\subsection{The sublattice $E_8(-1)^2$}\label{e8s}
It is not hard to identify a
sublattice of $NS(X_t)$ which is isomorphic to  $E_8(-1)^2$.
One copy of $E_8(-1)$ is spanned by the following eight $-2$-curves on $X_t$:
$$
\begin{array}{ccccccccccccc}
(x=0)&--&E_{001}^{x=0}&--&E_{001}&--&E_{001}^{z=0}&--&(z=0)&--&
E_{101}&--&(x=-1)\\
&&&&\mid&&&&&&&&\\
&&&&E_{001}^{x=-zt}&&&&&&&&
\end{array}
$$
another copy of $E_8(-1)$, perpendicular to this one, is given by:
$$
\begin{array}{ccccccccccccc}
E_{010}&--&E_{010}^{l_\infty}&--&(l_\infty)&--&E_{100}^{l_\infty}&--
&E_{100}&--&E_{100}^{z=-1}&--&(z=-1)\\
&&&&\mid&&&&&&&&\\
&&&&E_{110}&&&&&&&&
\end{array}
$$
By considering the effect of the symmetry on
the nodal curves one is led to the following expression for the
Nikulin involution.

\subsection{The Nikulin involution.}\label{2.5}
The pair of $E_8$'s in $NS(X_t)$ described in \ref{e8s} defines
a Nikulin involution $\iota$ on $X_t$
as in \cite{morrison}.
It is given by:
$$
\iota=\iota_t:X_t\longrightarrow X_t,\qquad
\iota(x,z,y)=\left(1/z,1/x,-y/(x^2z^2)\right)
$$
(the minus sign assures that $\iota$ has only isolated fixed points).

The invariants under the action of $\iota$
in the function field of the surface $X_t$,
are generated by
$$
\xi_1=x/z,\quad\xi_2=x+1/z\quad\mbox{ and }\quad\eta=y(xz-1)z^{-3}.
$$
The desingularization of $X_t/\iota$ is a $K3$ surface denoted by $V_t$:
$$
V_t:\qquad \eta^2=\xi_1(\xi_1+t)(\xi_1+\xi_2+1)(\xi_2^2-4\xi_1).
$$
Note that $\xi_1(\xi_1+t)(\xi_1+\xi_2+1)(\xi_2^2-4\xi_1)$ pulls back to
$xz(x+1)(z+1)(x+zt)$ times $\left((xz-1)/z^3\right)^2$.
As the isogeny is defined over $\QQ(t)$ we obtain:

\subsection{Lemma.}\label{vt}
{\sl
The desingularisation of the surface $X_t/\iota$ is the $K3$ surface $V_t$.
The graph of the rational map $X_t\rightarrow V_t$ defines
a correspondence, defined over $\QQ(t)$, which induces an isomorphism
on the transcendental parts of $H^2_\et$.
}

\subsection{An alternative description of the Nikulin involution.}
Using an elliptic fibration on $X_t$ given in \cite{AOP},
one obtains the following way to describe the involution $\iota$.

Consider the map
$$ \pi:\;X_t\cdots\rightarrow \PP^1\quad
(x,y,z)\mapsto \alpha:=\frac{y}{z(x+tz)}.
$$
This map in fact defines a morphism. The fibre over a general point
$\alpha\in{\PP}^1$ is the genus $1$ curve $D_{\alpha}$ with equation
$$
  \alpha^2z(x+zt)=x(x+1)(z+1).
$$
Using the change of coordinates
$$
\xi:=t\alpha^2/x,\quad\quad \eta:=(\xi+t\alpha^2)/z
$$
one obtains for $D_{\alpha}$ the equation
$$
\eta^2+(1-\alpha^2)\xi\eta+t\alpha^2\eta=\xi^3+t\alpha^2\xi^2.
$$
In this way, $\pi:\;X_t\to\PP^1$ is
the elliptic surface $\pi:\;D_{\alpha}\to\PP^1$ corresponding
to $(\xi,\eta,\alpha)\mapsto\alpha$.
Note that the fibre of this surface over $\alpha$ is
the same as the fibre over $-\alpha$.

Let $P_1$ be the section of this surface over $\PP^1$ given
by $P_1(\alpha)=(0,0,\alpha)$.
The Nikulin involution $\iota$ is then described as
$$
\iota(\xi,\eta,\alpha)=\left(P_1(\alpha)-(\xi,\eta),-\alpha\right),
$$
where $P_1(\alpha)-(\xi,\eta)$ is interpreted in terms of the group
law on the
elliptic curve $D_{\alpha}$.

\subsection{The branch locus of $V_t$.}
The branch locus of $V_t$ consists of four lines (including the line
at infinity) and a conic. The line $\xi_1+t=0$ meets the conic
transversely in two points, conjugate over the field $\QQ(\sqrt{-t})$,
whereas the other three lines are tangent to the conic and all
contain the (triple) point $(0,1,0)$. Blowing up the
singular points of the branch curve (in a point of tangency one must blow
up twice, in the triple point four times (see \ref{singbra})), one obtains
a rational surface such that $V_t$ is the double cover of this surface
branched over six disjoint smooth rational curves (the strict transforms
of the five irreducible components of the branch curve and the rational
curve which maps to the triple point). In the next section we will see
that one can blow down this rational surface to $\PP^2$ in such a way that
the images of these 6 rational curves are lines which are tangent to a conic.

\subsection{Remark.}
{}From \ref{Tgen} and \cite{morrison} we then have, for general $t$, that
$$
T_{V_t}\cong T(2)=\langle 4\rangle \oplus \langle 4\rangle \oplus \langle -4\rangle.
$$
This implies that the general $V_t$ is not isomorphic to the
Kummer surface of a product of two elliptic curves
(consider the transcendental lattices!).
It is not hard to check that for any elliptic curve
$E$ there is a subgroup $H\subset E\times E$, $H\cong (\ZZ/2\ZZ)^2$
such that $(E\times E)/H$ has transcendental lattice
$\langle 2\rangle \oplus \langle 2\rangle \oplus \langle -2\rangle=T$, hence the transcendental
lattice of the Kummer variety of $(E\times E)/H$ is $T(2)$. We will not use this
result explicitly since it does not guarantee the existence of a
correspondence over $\QQ(t)$.

\subsection{Five fold symmetry for $t=-1$.}
It is amusing to observe that in the case $t=-1$ one finds $25$ nodal curves
on the $K3$ surface $X_t$ which form a configuration already described by
Vinberg.

In case $t=-1$, the $6$ lines in $\PP^2$,
the $3$ exceptional divisors over the
double points ($(1:1:0)$, $(1:0:1)$ and $(0:1:1)$) and the $4\cdot 4=16$
curves over the $4$ triple points give a configuration of $25$ $-2$-curves
on $X_{-1}$.
The graph of this configuration (vertices correspond to the nodal curves,
edges are between vertices for which the corresponding nodal curves
intersect) is given in \cite{vinkap}, p.195, figure 2,
it has an obvious 5-fold symmetry!
The vertices in that figure are numbered from $1$ to $27$ with exception of
the numbers $19$ and $24$, the corresponding nodal curves can be chosen as:
$1=E_{100}$, $4=E_{111}^{z=-1}$, $5=E_{111}$, $6=E_{111}^{x=-1}$, $9=(z=0)$,
$11=E_{001}$, $15=(l_\infty)$, $17=E_{010}$, $23=(x=0)$,
from this it is easy to find the curves corresponding to the other vertices.

Note that the two copies of $E_8$ given in \ref{e8s} are exchanged by the
symmetry of the graph given by reflection in the vertical axis
(which contains the vertices 5, 27, 13 and 18).
A similar symmetry exists for general $t$
and is induced by the Nikulin involution.

\section{$V_t$ as a Kummer surface}\label{vkum}

\subsection{Lemma}\label{2.7}\label{4.1}
{\sl Let $W_t$ be the $K3$ surface defined by:
$$
W_t:\quad
t(t+1)\eta_2^2=(\xi_7^2+t\xi_8^2)(4\xi_7-4t\xi_8-t-1)
\left((\xi_7-2t\xi_8-t)^2+t(\xi_8+1)^2\right).
$$
There exists an isomorphism
$$
\phi:V_t\;\stackrel{\cong}{\longrightarrow}\; W_t
$$
which is defined over $\QQ(t,\sqrt{-t})$. Let
$\phi':V_t\rightarrow W_t$ be the
${\rm Gal}(\QQ(\sqrt{-t})/\QQ(t))$-conjugate of $\phi$,
and let $\Gamma_\phi$, $\Gamma_{\phi'}$ be their graphs in $V_t\times W_t$.
Then the correspondence $\Gamma_\phi+\Gamma_{\phi'}$, which is defined
over $\QQ$, induces an
isomorphism between the part of $H^2_\et(V_t)$ orthogonal to the
$19$ algebraically independent cycle classes and
the corresponding part of $H^2_\et(W_t)$.
}

\ts
To prove this, regard $V_{t}$ as a double cover of the plane,
with (affine) equation
$$\eta^2=\xi_1(\xi_1+t)(\xi_1+\xi_2+1)(\xi_2^2-4\xi_1).$$
 We will
explicitly
describe two Cremona transformations of the plane whose composition induces
the desired isomorphism $\phi$.

The ramification locus consists of $4$ lines (including the
line at infinity) and a conic;
note that $3$ of these lines (the lines $\xi_1=0$, $\xi_1+\xi_2+1=0$
and the line at infinity) are tangent to the conic.

We first apply the Cremona transformation which blows up these three
points of tangency and blows down the three lines connecting them.
In explicit (affine) coordinates, this map can be described by
$$
(\xi_1,\xi_2)\longmapsto(\xi_3,\xi_4):=\left(\xi_1(\xi_2+2)/(\xi_2^2-4\xi_1),
(\xi_2+2\xi_1)/(\xi_2^2-4\xi_1)\right).
$$
It transforms the three lines tangent to the conic
and the conic itself into four lines,
the remaining line (given by $\xi_1+t=0$) is mapped onto a conic.
One computes (the factors below correspond to the equations of the
resulting lines and conic)
$$
\eta_1^2 = \xi_3\xi_4(\xi_3+\xi_4+1)(2\xi_3^2+2t\xi_4^2+\xi_3+t\xi_4),
$$
where $\eta_1:=\eta\xi_2(\xi_2+2)(2\xi_1+\xi_2)(\xi_2^2-4\xi_1)^{-3}$.

In the coordinates $\xi_3,\xi_4,\eta_1$, this surface is again described as
a double cover of the plane ramified over a conic and $4$ lines, one of which
is the line at infinity.
Two of the lines intersect in
the point $(\xi_3,\xi_4)=(0,0)$
which is on the conic (hence the configuration has one triple point),
the other intersection points are ordinary double points.

Next, apply the Cremona transformation whose base points are this triple
point $(0,0)$, the point $\left(-t/(t+1),-1/(t+1)\right)$ in the intersection
of the line $\xi_3+\xi_4+1=0$ and the conic, and a (nonrational) point
$(s,1,0)$ where the line at infinity and the conic intersect
(note that $s^2=-t$).
This transformation has the property
that each of the $5$ components of the branch locus
has a line as image.

Explicitly, this second transformation can be given as
$$
(\xi_3,\xi_4)\mapsto (\xi_5,\xi_6):=\left((s-1)
\frac{\displaystyle \xi_3^2-s^3\xi_4^2+(s^2-s)\xi_3\xi_4}{\displaystyle
\xi_3+s^2\xi_4},
(s^2+s)\frac{\displaystyle (s-s^2)\xi_4^2+(s-1)\xi_3\xi_4+s\xi_4}
{\displaystyle \xi_3+s^2\xi_4}\right).
$$

It lifts to a birational map from our surface to the one given by
$$\begin{array}{c}t(t+1)\eta_2^2=\\
\xi_5\xi_6\left((1+s)\xi_5-s\xi_6+s^2+s\right)
\left((1-s)\xi_6+s\xi_5+s^2-s\right)
\left((2+2s)\xi_5+(2-2s)\xi_6+s^2-1\right),\end{array} $$
with
$\eta_2=\eta_1(1-s)(\xi_3-s\xi_4)\left((st+s)\xi_4-(t+1)\xi_3-t+s\right)
(\xi_3-t\xi_4)^{-2}$.

Finally, put
$$
\xi_7:=(\xi_5+\xi_6)/2,\quad \xi_8:=(\xi_5-\xi_6)/2s\qquad{\rm so}\quad
\xi_5=\xi_7+s\xi_8,\quad \xi_6=\xi_7-s\xi_8.
$$
With these coordinates, the equation is
$$
t(t+1)\eta_2^2=(\xi_7^2+t\xi_8^2)(4\xi_7-4t\xi_8-t-1)
\left((\xi_7-2t\xi_8-t)^2+t(\xi_8+1)^2\right),
$$
thus it defines a $K3$ surface, $W_t$,  which is defined over $\QQ(t)$.
The composition of the birational maps described so far yields
the isomorphism $\phi: V_t\;{\stackrel{\cong}{\longrightarrow}}\;W_t$,
defined over $\QQ(\sqrt{-t})$. Let $\phi'$ be the conjugate isomorphism
(defined by the same formulas as $\phi$ but with $-s$ for $s$).
A generator of $H^{2,0}(W_t)$ is given in local coordinates by
the regular $2$-form $\omega_W:=\dd\xi_7\wedge\dd\xi_8/\eta_2$.
A direct calculation shows that
$$
\phi^*\omega_W+(\phi')^*\omega_W\;=\;\dd\xi_1\wedge\dd\xi_2/\eta\neq 0.
$$
Hence the correspondence on $V_t\times W_t$, defined over $\QQ$,
which is the sum of the graphs $\Gamma_\phi+\Gamma_{\phi'}$
defines a nonzero map $H^{2,0}(W_t)\rightarrow H^{2,0}(V_t)$.
Thus it must induce an isomorphism on the transcendental lattices
of $W_t$ and $V_t$.
The comparison theorem for complex and $\ell$-adic cohomology
implies that the same is true for the corresponding Galois
representations.

This proves the lemma.\hfill{$\Box$}

\subsection{The $K3$ surface $W_t$.}\label{4.2}
The branch curve of the double cover $W_t\rightarrow\PP^2$
as described in Lemma~\ref{2.7} consists of $6$ lines (defined over
$\QQ(t,\sqrt{-t})$), including the line at infinity.
The smooth conic defined by
$4\xi_7+4t\xi_8^2-1=0$
is tangent to each of these lines. In particular, $W_t$ is the
Kummer surface of the Jacobian of the genus two curve $C_t$ which is
the double cover of the conic branched over the $6$ points of tangency
with the lines (see \cite[Exc.~VIII.6]{Beau}, \cite[\S3.10]{C-F}).
We briefly recall some of these classical results.

\subsection{Kummer surfaces and genus $2$ curves.}\label{Kmg2}
Let $K$ be a field of characteristic $\neq 2$. Suppose $C/K$ is
given by $y^2=f(x)$ for some separable polynomial $f\in K[x]$ of
degree $5$ or $6$. Over some extension field of $K$
we write $f(x)=\prod(x-a_j)$.

The Jacobian $JC$ of $C$ is birational to the symmetric product
$(C\times C)/S_2$, hence its function field is the subfield of
$K(x_1,x_2,y_1,y_2)$ (with the relations $y_i^2=f(x_i)$) of elements
fixed under the involution $\sigma$ given by $\sigma(x_1)=x_2$
and $\sigma(y_1)=y_2$. The Kummer surface $\Km(JC)$ is birational to
the quotient of $JC$ under the $[-1]$-map, hence its function field
is the subfield of $K(x_1,x_2,y_1,y_2)$ of elements fixed under
the two involutions $\sigma$ and $\iota$, with $\iota(x_i)=x_i$
and $\iota(y_i)=-y_i$.

The latter subfield is generated over $K$ by the functions
$\eta:=y_1y_2$ and $\xi:=x_1x_2$ and $\zeta:=x_1+x_2$.
They satisfy a relation
$$
\eta^2=F(\xi,\zeta)
$$
with $F$ the unique polynomial such that $f(x_1)f(x_2)=F(x_1x_2,x_1+x_2)$.

Observe that over an extension of $K$ one has
$$
f(x_1)f(x_2)=\prod((x_1-a_j)(x_2-a_j))=\prod(\xi-a_j\zeta+a_j^2).
$$
Hence one concludes that $\Km(JC)$ is birational over $K$ to a
double cover of the plane, ramified over six lines (including the
line at infinity in the case that the degree of $f$ is $5$).
The points $P_j:=(\xi=a_j^2,\zeta=2a_j)$ correspond to the pairs
of Weierstrass points $(T,T)\in C\times C$, with $T=(a_j,0)\in C$.
Note that the $P_j$ are also on the conic defined by $\zeta^2=4\xi$,
and in fact the line $\xi-a_j\zeta+a_j^2=0$ is tangent to this
conic in $P_j$. The same is true for the line at infinity (the point
of tangency comes from the point at infinity on $C$ in the case
where the degree of $f$ is $5$).
This shows that seen as a double cover of the plane, $\Km(JC)$
is ramified over six lines which are tangent to a given conic.
The inverse image of this conic is given by the two equations
$\zeta^2=4\xi$ and $\eta^2=\prod(\zeta/2-a_j)^2$. Hence it consists
of two irreducible components, both defined over $K$.
Moreover, we can recover the Weierstrass points of $C$ (and hence
$C$ itself up to a quadratic twist) from the six points of tangency of
the lines with the conic.

\subsection{Lemma.}\label{kjct}{\sl
The $K3$ surface $W_t$ studied in \ref{4.1} and \ref{4.2} is isomorphic
to $\Km(JC_t)$, where the genus two curve $C_t$  is defined by
$$
C_t:\quad y^2=x(x^2-4x+4+4t)(x^2+4x+4+4t).
$$
This isomorphism is defined over $\QQ(t,\sqrt{-t}, \sqrt{-t-1})$.

Given an equation $\eta^2=F(\xi,\zeta)$ for $\Km(JC_t)$ as above,
let $\Km(JC_t)^{(-t-1)}$ be the `twist'
defined by $(-t-1)\eta^2=F(\xi,\zeta)$.

Then there is a correspondence on $W_t\times \Km(JC_t)^{(-t-1)}$,
defined over $\QQ(t)$, which induces an isomorphism of
$G_{\QQ(t)}$-representations between the transcendental
parts of the $H^2_{\et}$'s.
}

\ts
As before, write $s^2=-t$. We will use new coordinates to describe $W_t$,
namely $\xi_9$ and $\xi_{10}$ given by
$$
\xi_8=\frac{\xi_{10}}{8s}-\frac12
$$
and
$$
\xi_7=(\xi_9-2s\xi_{10}-4t+4)/16.
$$
In these coordinates, the conic $4\xi_7+4t\xi_8^2-1=0$ becomes
$\xi_{10}^2=4\xi_9$.
The $6$ lines over which $W_t\to {\mathbb P}^2$ ramifies become
the line at infinity and five lines $\xi_9-b_j\xi_{10}+b_j^2=0$,
with
$$
\{b_1,b_2,b_3,b_4,b_5\}\;=\;\{0, 2+2s, 2-2s,-2+2s,-2-2s\}.
$$
The equation for $W_t$ in the new coordinates is
$$
W_t:\quad 2^{18}t(t+1)\eta_2^2=\prod_{j=1}^5 (\xi_9-b_j\xi_{10}+b_j^2),
$$
which is in fact an equation over $\QQ(t)$.

The discussion in \ref{Kmg2} above shows that provided we have a
square root of $t(t+1)$ available, this defines a birational model of
the Kummer surface $\Km(JC_t)$ where $C_t$ is the hyperelliptic curve
with Weierstrass points over infinity and over the $b_j$'s, so the equation
defining $C_t$ is the one given in the lemma:
$$
y^2=\prod_{j=1}^5 (x-b_j)=x(x^2-4x+4+4t)(x^2+4x+4+4t).
$$

To show the second part, put
$$
\Km(JC_t)^{(-t-1)}:\quad (-t-1)\eta^2=\prod_{j=1}^5 (\xi-b_j\zeta+b_j^2).
$$
A birational map $\psi:\Km(JC_t)^{(-t-1)}\to W_t$ is given by
$$
\psi(\eta,\xi,\zeta):=\left(\eta_2=\frac{2^{-9}\eta}{s},
\xi_7=\frac{(\xi-2s\zeta-4t+4)}{16},\xi_8=\frac{\zeta}{8s}-\frac12\right).
$$
One obtains the `conjugate' $\psi'$ by replacing all occurrences of
$s$ by $-s$ in this description. A direct calculation reveals that
$$
\psi^*\frac{\dd \xi_7\wedge \dd \xi_8}{\eta_2}=
\psi'^*\frac{\dd \xi_7\wedge \dd \xi_8}{\eta_2}=
4\frac{\dd \xi\wedge \dd \zeta}{\eta},
$$
from which the lemma follows by the same argument as in the proof
of Lemma~\ref{4.1}.
\qed

\subsection{The product of elliptic curves.}\label{4.5}
The curve $C_t$ has, besides the hyperelliptic involution,
another involution:
$$
\varphi=\varphi_t:C_t\longrightarrow C_t,\qquad
\varphi(x,y):=(r^2/x,r^3y/x^3)\qquad(r^2=4+4t).
$$
The quotient by this involution is an elliptic curve.
In fact, the invariant functions on $C_t$ are generated by
$\eta:=y(x+r)/x^2$ and $\xi:=-x/(2r)-r/(2x)$
and the quotient curve $F_t$ is defined by
$$
F_t:=C_t/\varphi:\quad \eta^2=-8r^3(\xi-1)(\xi^2-\tpi).
$$
Replacing $r$ by $-r$ yields yet another involution (namely, the
composition of the previous one and the hyperelliptic involution $\tau$)
and hence a second elliptic curve
$$
F'_t:=C_t/(\varphi\circ\tau):\quad \eta^2=8r^3(\xi-1)(\xi^2-\tpi).
$$

By considering the pull back to $C_t$ of the invariant
differentials on these elliptic curves one concludes that
$JC_t$ is isogenous (over $\QQ(t,r)=\QQ(t,\sqrt{t+1})$)
to the product $F_t\times F'_t$ of these two elliptic curves.
In explicit form, this isogeny is obtained from the
two quotient maps $\alpha:C_t\to F_t$ and $\alpha':C_t\to F_t'$
using
$$
C_t\times C_t\longrightarrow F_t\times F_t'\;\quad
(P,Q)\longmapsto \left(\alpha(P)+\alpha(Q),\alpha'(P)+\alpha'(Q)\right).
$$

The associated Kummer surfaces are isogenous (again, over $\QQ(t,r)$
and not necessarily over $\QQ(t)$!) as well.
It is easily seen that the Kummer surface of $F_t\times F_t'$
is birational over $\QQ(t,r)$ to the surface
with equation
$$
   -y^2=(x_1-1)(x_1^2-\tpi)(x_2-1)(x_2^2-\tpi).
$$
By twisting, this also gives a rational map, defined over $\QQ(t,r)$,
 from $\Km(JC_t)^{(-t-1)}$
to the surface defined by
$ -(-t-1)y^2=(x_1-1)(x_1^2-\tpi)(x_2-1)(x_2^2-\tpi)$.
Note that the latter equation in fact defines
the Kummer surface of $E_t\times E_t^{(t+1)}$ over $\QQ(t)$.

Now we show that this rational map together
with its $\QQ(t,r)/\QQ(t)$-conjugate yields a correspondence
defined over $\QQ(t)$ between $\Km(JC_t)^{(-t-1)}$
and $\Km(E_t\times E_t^{(t+1)})$ with the property that
it is nonzero on the transcendental part of $H^2_{\et}$'s.

\subsection{Genus $2$ curves with non-simple Jacobians}
Suppose $k$ is a field of characteristic $\neq 2$. Let $f\in k[x]$
be a separable polynomial of degree $5$  and
$C:\;y^2=f(x)$. The regular differentials on $C$ form a
$k$-vector space with a basis $\frac{\dd x}{y},\;x\frac{\dd x}{y}$.
Assume that for $i=1,2$ an elliptic curve $E_i$ over $k$ is given,
with a nonzero regular differential $\omega_i$ on $E_i$
and a morphism
$$\alpha_i:\;C\longrightarrow E_i$$
having the property that $\alpha_1^*\omega_1$ and $\alpha_2^*\omega_2$
are linearly independent. Moreover we assume that $\alpha_i$
sends the point at infinity on $C$ to the zero on $E_i$. Write
$\alpha_i^*\omega_i=(a_ix+b_i)\frac{\dd x}{y}$.
The independence of the pull backs can be phrased by saying that
$$
d:=a_1b_2-a_2b_1\neq 0.
$$

Consider the commutative diagram of rational maps
$$ \begin{array}{ccccc}
  C\times C && \cdots\cdots{\longrightarrow} &&\Km(JC) \\
   {\Big\downarrow}\psi &&&& {\Big\downarrow}\\
  E_1\times E_2 && \cdots\cdots{\longrightarrow} &&\Km(E_1\times E_2).
\end{array}
$$
Here $\psi$ is the morphism $\psi:\;(P,Q)\mapsto\left(
\alpha_1(P)+\alpha_1(Q),\alpha_2(P)+\alpha_2(Q)\right)$.

Note that $\omega_1\wedge\omega_2$ can be regarded both as a
regular $2$-form on $\Km(E_1\times E_2)$ and as the regular $2$-form
on $E_1\times E_2$ obtained as the pull back of the one on the Kummer.
One computes that
$$
\psi^*(\omega_1\wedge\omega_2)=d(x_1-x_2)\frac{\dd x_1\wedge\dd x_2}{y_1y_2}
$$
using coordinates $x_1,y_1,x_2,y_2$ on $C\times C$ which satisfy
$y_i^2=f(x_i)$.

As is explained in (\ref{Kmg2}) above, one has
coordinates $\eta=y_1y_2$ and $\xi=x_1x_2$ and $\zeta=x_1+x_2$
on $\Km(JC)$. The regular $2$-form
$\eta^{-1}\dd\xi\wedge\dd\zeta$
on $\Km(JC)$ pulls back under the horizontal rational map at
the top of the diagram above to
$
(y_1y_2)^{-1}\dd(x_1x_2)\wedge\dd(x_1+x_2)=
(x_2-x_1)(y_1y_2)^{-1}\dd x_1\wedge\dd x_2
$.

Combining the above pull backs, one concludes that using the
vertical arrow on the right of our diagram, $\omega_1\wedge\omega_2$
pulls back to
$-d\eta^{-1}\dd\xi\wedge\dd\zeta$ on $\Km(JC)$.

\subsection{} \label{4.7}
We now apply this to the situation described in (\ref{4.5}).
Here we have
$$\rho:\; \Km(JC_t)\longrightarrow \Km(F_t\times F_t')\cong V_t$$
where $V_t$ is defined by $-y^2=(x_1-1)(x_1^2-\tpi)(x_2-1)(x_2^2-\tpi)$.
The surface $\Km(F_t\times F_t')$ corresponds to the
equation $-\tilde{y}^2=64r^6(x_1-1)(x_1^2-\tpi)(x_2-1)(x_2^2-\tpi)$.
An isomorphism between $V_t$ and this surface is described by
$\tilde{y}=8r^3y$. Hence the $2$-form
$y^{-1}\dd x_1\wedge\dd x_2$
on $V_t$ pulls back to
$(8r^3\tilde{y})^{-1}\dd x_1\wedge\dd x_2$
on $\Km(F_t\times F_t')$.

Since $\alpha^*\frac{\dd \xi}{\eta}=
\left(\frac{rx}{8t+8}+\frac12\right)\frac{\dd x}{y}$
and $\alpha'^*\frac{\dd \xi}{\eta}=
\left(\frac{-rx}{8t+8}+\frac12\right)\frac{\dd x}{y}$,
it follows that
$$ d=\frac{r}{8t+8}. $$
Hence, the pull back of $\frac{\dd x_1\wedge\dd x_2}{8r^3\tilde{y}}$
to $\Km(JC_t)$ is $-\frac{\dd\xi\wedge\dd\zeta}{(16t+16)^2\eta}$.

One concludes that
$$
\rho^*\left(\frac{\dd x_1\wedge\dd x_2}{y}\right)=
-\frac{\dd\xi\wedge\dd\zeta}{(16t+16)^2\eta}.
$$
Denoting by $\rho'$ the $\mbox{Gal}(\QQ(t,r)/\QQ(t))$-conjugate of $\rho$,
it follows as before that the sum of the graphs
$\Gamma_{\rho}+\Gamma_{\rho'}$ defines a correspondence
over $\QQ(t)$ which is nonzero on the transcendental parts of $H^2_{\et}$'s.

Twisting all surfaces over $\QQ(t,\sqrt{-t-1})$ one obtains the
same conclusion for the surfaces $\Km(JC_t)^{(-t-1)}$
and $\Km(E_t\times E_t^{(t+1)})$, hence we proved:

\subsection{Lemma.}
There is a correspondence on
$\Km(JC_t)^{(-t-1)}\times Km(E_t\times E_t^{(t+1)})$,
defined over $\QQ(t)$, which induces an isomorphism of
$G_{\QQ(t)}$-representations between the transcendental
parts of the $H^2_{\et}$'s.

\subsection{Conclusion.}
Putting together the various correspondences we have constructed,
one obtains
the desired correspondence defined over $\QQ(t)$ on the product of $X_t$ and
$\Km(E_t\times E^{(t+1)}_t)$. This finishes the proof of Theorem~\ref{main}.

The maps we constructed, over finite extensions of $\QQ(t)$,
compose (eventually after a further field extension to undo twists)
to give a dominant rational map:
$$
X_t\longrightarrow V_t\longrightarrow W_t\longrightarrow
K(JC_t)\longrightarrow K(E_t\times E_t),
$$
hence also Proposition~\ref{prp} follows.

\section{The fibre at $t=-1$} \label{a1}

\subsection{}
We conclude this paper with some remarks on various models of the famous
$K3$ surface $X_{-1}$. As observed in \ref{x1},
the $K3$ surface $X_{-1}$
is the desingularisation of the quotient of
$E_i\times E_i$, the self product of the elliptic curve $E_i=\CC/\ZZ[i]$,
by the automorphism $\phi$ of order $4$ induced by
$(z_1,z_2)\mapsto (iz_1,-iz_2)$ on $\CC^2$, \cite{SI}.
Below we show how to obtain this isomorphism directly from the equations
defining $X_{-1}$ and $E_i$. It is convenient to use
projective coordinates $(u:v:w)=(x:z:-1)$, so the equation for
$X_{-1}$ becomes:
$$
X_{-1}:\qquad \sigma^2=uvw(u-w)(v-w)(u-v).
$$

\subsection{}
The elliptic curve $E_i$ is isomorphic to $E:\;t^2=s(s^2-1)$ and also to
$E':\;y^2=x(1-x^2)$, hence $E_i\times E_i \cong E\times E'$ and $\phi$
may be given by $\phi((s,t),(x,y))=((-s,it),(-x,-iy))$.
The quotient map $E\times E'\rightarrow X_{-1}$ is given by
$$
(\sigma:u:v:w)=(xs(xs+1)(x+s)ty: \; xs^2-x:\,xs^2+s:\,xs^2+s^2x)
$$
in fact, a direct calculation shows that
$$
uvw(u-v)(u-w)(v-w)=(xs(xs+1)(x+s))^2x(1-x^2)s(s^2-1).
$$
This map was found from the results below.

\subsection{Vinberg's model.}
The surface $X_{-1}$ has a projective model $Y$ which is
a singular quartic surface
in $\PP^3$ (see \cite{vinberg}, Theorem 2.5, but we replaced
$X_0$ there by $\zeta X_0$ for a $\zeta\in k$ with $\zeta^4=-1$):
$$
Y:\quad X_0^4=X_1X_2X_3(X_1+X_2+X_3).
$$
The elliptic curve $E_i$ is isomorphic to $E:\;t^2=s^4-1$ and also to
$E':\;y^2=x^4+1$, hence $E_i\times E_i\cong E\times E'$ and $\phi$
may be given by $\phi((s,t),(x,y))=((is,t),(-ix,y))$.
The quotient map $E\times E'\rightarrow Y$ is given by
$$
(X_0:X_1:X_2:X_3)=(sx:y-1:1+t:1-t),
$$
it is easy to see that this map has degree $4$ and is invariant under $\phi$.
This map was found by studying the pencil of curves on $Y$ defined by
$X_3=\lambda X_2$.

\subsection{}
An isomorphism $X_{-1}\rightarrow Y$ is given by
$$
(X_0:X_1:X_2:X_3)=(\sigma:vw(v-w):-uw(u-w):uv(u-v)).
$$
Note that  $X_1+X_2+X_3=(u-v)(u-w)(v-w)$ and thus the equation for $Y$
pulls back to $\sigma^4=(uvw(u-v)(u-w)(v-w))^2$.

\

\end{document}